\documentstyle{amltd2004}
\begin{document}
\annalsline{155}{2002}
\received{September 6, 2000}
\startingpage{189}
\def\bye{\end{document}}
 \font\tenrm=cmr10
\def\thefootnote{\fnsymbol{footnote}}
\input amssym.tex
\input amssym.def

 \def\bcw{\mathbin{\bigcirc\mkern-15mu\wedge}}

 \def\r#1{{\mathop{#1}\limits^\circ}}
\def\ring{\r}

 \def\nint{\mathbin{\int\mkern-18mu\diagup \;} }

\newcommand{\namelistlabel}[1] {\mbox{#1}\hfil}
\newenvironment{namelist}[1]{%
\begin{list}{}
{\let\makelabel\namelistlabel
\settowidth{\labelwidth}{#1}
\setlength{\leftmargin}{1.1\labelwidth}}
}{%
\end{list}}
%--------------- Author macros ---------------
%for Bbb in amstex
\catcode`\@=11
\font\twelvemsb=msbm10 scaled 1100
\font\tenmsb=msbm10
%\font\ninemsb=msbm7 scaled 1100%msbm9
\font\ninemsb=msbm10 scaled 800
\newfam\msbfam
\textfont\msbfam=\twelvemsb  \scriptfont\msbfam=\ninemsb
  \scriptscriptfont\msbfam=\ninemsb
\def\msb@{\hexnumber@\msbfam}
\def\Bbb{\relax\ifmmode\let\next\Bbb@\else
 \def\next{\errmessage{Use \string\Bbb\space only in math
mode}}\fi\next}
\def\Bbb@#1{{\Bbb@@{#1}}}
\def\Bbb@@#1{\fam\msbfam#1}
\catcode`\@=12

 \catcode`\@=11
\font\twelveeuf=eufm10 scaled 1100
\font\teneuf=eufm10
\font\nineeuf=eufm7 scaled 1100%eufm9
\newfam\euffam
\textfont\euffam=\twelveeuf  \scriptfont\euffam=\teneuf
  \scriptscriptfont\euffam=\nineeuf
\def\euf@{\hexnumber@\euffam}
\def\frak{\relax\ifmmode\let\next\frak@\else
 \def\next{\errmessage{Use \string\frak\space only in math
mode}}\fi\next}
\def\frak@#1{{\frak@@{#1}}}
\def\frak@@#1{\fam\euffam#1}
\catcode`\@=12
%-------------- Author entries --------------------

\title{Discrete analogues in harmonic analysis:\\ Spherical averages} %Article title
\shorttitle{Discrete analogues in harmonic analysis} % Shortened version for headline title

  \acknowledgements{This work was supported by NSF Grants: DMS-9970899, DMS-9706889, and 
DMS-9731647.}
  \twoauthors{A. Magyar, E. M. Stein,}{S. Wainger}
 \institutions{University of Wisconsin-Madison, Madison, 
WI\\
{\eightpoint {\it E-mail addresses\/}: magyar@math.wisc.edu}\\
\hglue.97in {\eightpoint wainger@math.wisc.edu} \vglue6pt
Princeton University, Princeton, NJ\\
{\eightpoint {\it E-mail address\/}: stein@math.princeton.edu}}
%-------------- Article Text--------------------

\bigbreak \centerline{\bf Abstract}
\vglue9pt 
In this paper we prove an analogue in the discrete setting of ${\Bbb{Z}^d}$, of
the spherical maximal theorem for ${\Bbb{R}^d}$.  The methods used are two-fold:
the application of certain ``sampling'' techniques, and ideas arising in the 
study of the number of representations of an integer as a sum of $d$ squares,
in particular, the ``circle method''.  The results we obtained are by
necessity limited to $d \geq 5$, and moreover the range of $p$ for the $L^p$
estimates differs from its analogue in~${\Bbb{R}^d}$.
\vglue12pt

\section{Introduction}

 Geometric considerations, in particular curvature, play an
important role in harmonic analysis in ${\Bbb{R}^d}$.  Emblematic
of this are the properties of the spherical maximal function. 
Given the significance of this operator, it is an interesting
and natural question to ask what happens when we consider its
discrete analogue; that is, what can be said of the corresponding
version of the spherical maximal theorem taken over ${\Bbb{Z}^d}$?
It is the purpose of this paper to answer this question by proving
optimal $\ell^p$ estimates in this setting.

We shall now describe these results, turning first to
${\Bbb{R}^d}$. The spherical averages are defined by the operators
${\cal{A}}_\lambda$,where
$$
{\cal{A}}_\lambda \, ( f ) \, = \, f \, \star \, d \sigma_\lambda $$
with $d \sigma_\lambda$ the normalized invariant measure on the 
sphere $| x | = \lambda$.  With the definition of the
maximal function, ${\cal{A}_\star} ( f ) ( x ) =
\displaystyle{\sup_{0 < \lambda < \infty}} | 
{\cal{A}_\lambda} ( f ) ( x ) |$, we recall the main estimate for it, 
\begin{equation}
\left\| {\cal{A}}_\star \, ( f ) \, \right\|_{L^p ( {\Bbb{R}^d})} \, \leq \,
A \, \left\| f \right\|_{L^p ( {\Bbb{R}^d})} \, ,
\ \ {\rm if } \ \ p > \frac{d}{d - 1} \ \ {\rm and } \ \ d \, \geq 2 \, . \hskip.5in
\end{equation} 
(See [S], [SW1], [B1].)
 
The discrete analogue of ${\cal{A}}_\lambda$ is the operator
\begin{equation}
A_\lambda \, ( f )( n ) \, = \, \frac{1}{N ( \lambda )} \;
\displaystyle{\sum_{|m| = \lambda}} \; f ( n - m ) \, .
\end{equation}
Here $n$ and $m$ are restricted to range over ${\Bbb{Z}^d}$; also
$N ( \lambda ) =$ the number of $m \in {\Bbb{Z}^d}$, so that 
$| m | = \lambda$.  Notice that only those $\lambda$ 
for which $\lambda^2$ is an integer are relevant; also
observe that $N ( \lambda ) = r_d ( \lambda^2)$, where
$r_d ( k )$ is the standard counting function giving the number
of ways of representing $k$ as a sum of $d$ squares.

Now, up to this point, formulating a discrete analogue of the
spherical maximal function, i.e.
$A_\star ( f ) ( n ) = 
\displaystyle{\sup_{0 < \lambda < \infty}} | A_\lambda ( f ) ( n )|$, and
asking the question of its $\ell^p$ boundedness, 
seem quite straightforward.\footnote  
{Here we use the notation that for $f$ defined on
${\Bbb{Z}^d}$, it belongs to 
$\ell^p ( {\Bbb{Z}^d}) $ if
$\displaystyle{\sum_{n \in {\Bbb{Z}^d}}} | f ( n ) |^p$ is finite.
The $\ell^p$ norm is of course the $p^{\rm th}$ root of the 
sum.}
However, this is misleading since quite different ideas must come into play in
the discrete analogue, and anyway, the range of exponents is not the same
as the version in ${\Bbb{R}^d}$.  The theorem we prove is the following
optimal result.

\nonumproclaim{Theorem}  The maximal operator $A_\star$ is bounded in $\ell^p ({\Bbb{Z}^d})$ to
itself for $p > \frac{d}{d - 2}${\rm ,} when $d \geq 5${\rm .}
\endproclaim

Alex Ionescu has pointed out to us that simple examples show that this
result cannot be improved: in fact, when $d \geq 5$, $A_\star$
is not bounded on $\ell^p$ for $p \leq \frac{d}{d - 2}$; moreover when
$ d < 5$, the $\ell^p$ boundedness fails for every $p < \infty$  (the 
case $p = \infty$ is of course trivial).  The relevant examples can be
found in Section~8.  Here the facts that the number of representations $r_d (k)$ 
is an irregular function of $k$ when $d \leq 4$, while 
$r_d ( k ) \approx k^{\frac{d-2}{2}}$ when $d \geq 5$, play  a role.
(For these assertions about $r_d ( k )$, consult [W].)
 
Our attack on the discrete spherical maximal function proceeds in
three stages.  To begin with (motivated by the ideas of the
circle method) we approximate $A_\lambda$ by an infinite sum of
simpler operators 
\begin{equation}
M_\lambda \, = \, c_d \, \sum \, 
e^{-2 \pi i \lambda^{2}a/q}  
\, M^{a/q}_\lambda \, ,
\end{equation}
with each $M^{a/q}_\lambda$ associated to a reduced fraction 
$a/q$, with $0 < a/q \leq 1$.  Now since each $M^{a/q}_\lambda$ is
a convolution operator on ${\Bbb{Z}^d}$, it 
corresponds to a Fourier multiplier $m^{p/q}_\lambda ( \xi )$, which
is given by
$$
m^{a/q}_\lambda ( \xi ) \, = \,
\displaystyle{\sum_{\ell \in {\Bbb{Z}^d}}} \,
G ( a / q , \ell ) \,
\Phi_q ( \xi \, - \ell/q ) \;
d \hat{\sigma}_\lambda ( \xi \, - \ell/q ) \, .
$$
Here $G$ is a normalized Gauss sum, $\Phi_q$ is a 
suitable cut-off function, and $d \hat{\sigma}_\lambda$ is the
Fourier transform of the unit measure $d \sigma_\lambda$ on the
sphere $| x | = \lambda$.

Notice that the first term of $M_\lambda$, 
corresponding to $a/q = 1 \equiv 0 \, {\rm mod} \, 1$, can be viewed as the
vestige of the continuous analogue on ${\Bbb{R}^d}$.  All the other
terms are approximations corresponding to the other rationals.

The second stage is to study each $M^{a/q}_\lambda$ as a sort
of discrete analogue of an operator on ${\Bbb{R}^d}$.  The main tool
is a general abstract theorem which allows one to pass from
certain convolution operators on ${\Bbb{R}^d}$ to 
analogous operators on ${\Bbb{Z}^d}$.  While ideas about 
special cases of this principle have been implicit in the past, 
our general
approach seems both new  interesting in its own right.  It is
presented in Section 2.  It is based in part on variants of
``sampling'' ideas which go back to Plancherel and P\'{o}lya [PP] and
which were taken up again later by Shannon [ShW].  Using arguments of a 
different kind,   Bourgain obtained certain results of this form; see  [B2, (3.5)].

The final stage of the argument is to show that $M_\lambda$ is an
adequate approximation of $A_\lambda$.  This is begun in
Sections 4 and 5, and is concluded in Sections 6 and 7.

The analysis of our theorem has as its starting point
a partial result obtained previously by one of us [M] (see
Proposition (4.2) below).  The interested reader may also want to 
compare the related ways the sums
$\displaystyle{\sum_{|n| = \lambda}} \, e^{2 \pi i n \cdot \xi}$
are treated in our paper (see Section 5), 
and in a previous work of Bleher and Bourgain [BB,  \S 6].  The context of that paper is however quite different
from ours.

\section{Discrete analogues of convolution operators}
 \advance\eqcount by -1

 Suppose $T ( f ) = f \star K$ is a convolution operator in 
$L^p (\Bbb{R}^d)$ to itself with a suitable distribution kernel $K$. 
Then, as is known, its Fourier transform 
$\widehat{K} = m ( \xi ) = 
\displaystyle{ \int_{\Bbb{R}^d}} 
\, K ( x ) e^{-2 \pi i x \xi} dx$ is a bounded function, 
and we can think of $T$ as a Fourier multiplier operator given by 
$( T f )^\wedge ( \xi ) = m ( \xi ) \hat{f} ( \xi)$.

To be precise, in what follows we shall assume in this section that in addition to $m ( \xi)$ being bounded,
  it is supported in the fundamental cube $Q = \{ \xi = ( \xi_j ) : - 1/2 < \xi_j \leq 1/2$,
$j = 1 , \ldots \ell \}$.  In this case $K ( x ) = \displaystyle{\int_{\Bbb{R}^d}} \,
e^{2 \pi i x \xi} \, m ( \xi ) d \xi$ is an $L^2$ 
function on ${\Bbb{R}}^d$, which is continuous (in fact, $C^\infty$).  
Thus $K_{\rm dis} = K \big|_{ \Bbb{Z}^d}$ is well-defined, as is
the convolution operator acting on functions on ${\Bbb{Z}^d}$ given by
$$
T_{\rm dis} ( f ) \, = \,
f \star K_{\rm dis} \, , \hspace{.20in}  T_{\rm dis} ( f ) ( n ) \, = \,
\displaystyle{\sum_{m \in {\Bbb{Z}}^d}} \,
K ( m ) \, f \, ( n - m ) \, . 
$$

Note that the condition that the multiplier be supported in $Q$ is natural.  Because then
not only does $T$ determine $T_{\rm dis}$, but conversely 
$T_{\rm dis}$ determines $T$, i.e. $K \big|_{\Bbb{Z}^d}$ determines $K$.
This follows since $K \big|_{\Bbb{Z}^d}$ determines the Fourier coefficients of the function 
$m ( \xi )$, when expanded as a function on $Q$.

Let $m_{\rm per}$ be the periodic extension of $m$, i.e. 
$m_{\rm per} ( \xi ) = \displaystyle{\sum_{\ell \in {\Bbb{Z}^d}}} \, 
m ( \xi - \ell )$.
Then $m_{\rm per} ( \xi )$ is the Fourier multiplier corresponding to $T_{\rm dis}$ in the sense
$$
\displaystyle{\sum_{n \in \Bbb{Z}^d}} \,
T_{\rm dis} \, ( f ) ( n ) \,
e^{- 2 \pi i n \xi} \, = \,
m_{\rm per} ( \xi ) \,
\displaystyle{\sum_{n \in \Bbb{Z}^d}} \,
f ( n ) e^{- 2 \pi i n \xi} \, ,
$$
for suitable functions $f$ on ${\Bbb{Z}^d}$.

Let us note that
\setcounter{equation}{-1}
\begin{equation}
m_{\rm per} ( \xi) \, = \,
\displaystyle{\sum_{n \in {\Bbb{Z}^d}}} \, K ( n ) \,
e^{- 2 \pi i n \xi} 
\end{equation}
in the sense of $L^2$ convergence of the series on any compact subset of ${\Bbb{R}^d}$.

In fact, $m ( \xi) = \displaystyle{\sum_{n \in {\Bbb{Z}^d}}} \,
K ( n ) \, e^{- 2 \pi i n \xi}$ on $Q$ 
represents the Fourier inversion of the
identity $K ( n ) = \int m ( \xi) e^{2 \pi i n \xi} \, d \xi$ 
$( m ( \xi )$ is supported in $Q$); and, moreover,  $m_{\rm per} ( \xi )$ is the periodic function which agrees with $m (
\xi )$ on $Q$.  This  establishes~(2.0).

The question we will be concerned 
with is how the norm of $T_{\rm dis}$ as an operator on
$\ell^p ({\Bbb{Z}^d})$ is controlled by the norm of the operator $T$ acting on $L^p ( \Bbb{R}^d)$. For our
applications it will be important to be able to deal with the more general case where the
$L^p$ and $\ell^p$ spaces of complex-valued functions are replaced by the spaces
$L^p_B ( {\Bbb{R}^d})$ and $\ell^p_B ( {\Bbb{Z}^d} )$ of functions taking their values in
the Banach space $B$.  In order to avoid technical problems involving definability, measurability, 
etc., we shall restrict our attention to the case when the 
Banach spaces in question are finite-dimensional.
However, all our estimates will be independent of the Banach spaces in question, so that a limiting argument will
encompass the results in the generality needed.  In particular, this argument will apply to the case when $B$ is an
$L^\infty$ space, which is what is needed for the maximal theorems below.

We shall suppose that $B_1$ and $B_2$ are a pair of finite-dimensional 
Banach spaces, and assume that 
$m ( \xi )$ is a bounded measurable function, taking its values in
${\cal{L} }( B_1 , B_2 )$; and as we have said we suppose $m$ is supported in $Q$.  Then $T$, described above, is a
bounded mapping from $L^2_{B_1} ( {\Bbb{R}^d})$ to $L^2_{B_2} ( {\Bbb{R}^d})$, and similarly $T_{\rm dis}$ is
bounded from $\ell^2_{B_1} ( {\Bbb{Z}^d})$ to $\ell^2_{B_2} ( {\Bbb{Z}^d})$.

\proclaim{Proposition}
Fix $p${\rm ,} $1 \leq p \leq \infty${\rm .}  If $T$ is bounded from 
$L^p_{B_1} ( {\Bbb{R}^d})$ to 
$L^p_{B_2} ( {\Bbb{R}^d} )${\rm ,} then $T_{\rm dis}$ is bounded from $\ell^p_{B_1} ( {\Bbb{Z}^d})$ to 
$\ell^p_{B_2} ( {\Bbb{Z}^d} )${\rm .  }
For these operators we have the norm inequality \pagebreak

\phantom{hi}
\vglue-28pt
\begin{equation}
\big\| \, T_{\rm dis} \, \big\|_{\ell^p_{B_1} \rightarrow \ell^p_{B_2} } \,
\leq \, 
C \big\| \, T \, \big\|_{L^p_{B_1} \rightarrow L^p_{B_2}} \, , 
\end{equation}
with a bound $C$ that depends only on the dimension $d${\rm ,} but not on $p$ or the Banach spaces $B_1$ and $B_2${\rm .}
\endproclaim

\demo{{R}emarks}
%\begin{namelist}{xxxx}
%\item[(1)] 
(1) It would be interesting to know if $C$ can be taken to be independent 
of the dimension $d$, or for that matter if $C = 1$.
%\item[(2)] 
\medbreak
(2) There is a converse to (2.1); i.e., a reverse inequality also holds.  
Since that fact will not be used below we will omit its
proof.
%\end{namelist} 
\enddemo

The proof of the proposition requires the following ``sampling'' and extension lemma.  We fix
the function $\Psi$ on ${\Bbb{R}^d}$ by  
$$
\Psi ( x ) \, = \,
\left(
\frac{\sin \pi x_1}{\pi x_1} \right)^2 \;
\left(
\frac{\sin \pi x_2}{\pi x_2} \right)^2 \, \ldots \,
\left(
\frac{\sin \pi x_d}{\pi x_d} \right)^2 \, ,\quad  x \, = \, ( x_1 , x_2 , \ldots x_d ) \, .
$$ 
For any suitable function $f$ on ${\Bbb{Z}^d}$ we consider its extension $f_{\rm ext} = F$ on ${\Bbb{R}^d}$ given by
\begin{equation}
F ( x ) \, = \,
f_{\rm ext} ( x ) \, = \,
\displaystyle{\sum_{n \in {\Bbb{Z}^d}}} \,
f ( n ) \, \Psi ( x - n ) \, .
\end{equation} 
(Note that if $f \in \ell^p$ for some $p$, the series above 
converges for every
$x \in {\Bbb{R}^d}.)$  We observe that in fact $F \big|_{\Bbb{Z}^d} = f$, since 
$\Psi ( 0 ) = 1$, and 
$\Psi ( n ) = 0$, if $n \in {\Bbb{Z}^d}$, $n \ne 0$; thus $f_{\rm ext}$ is a genuine
extension of $f$.  
The following estimate holds for any (finite-dimensional) Banach space $B$.

\specialnumber{2.1}\proclaim{Lemma}
If $f \in \ell^p ( {\Bbb{Z}^d} , B)${\rm ,} then $F \in L^p ( {\Bbb{R}^d} , B)${\rm ,} and
\begin{equation}
(1/A) \, \left\| f \right\|_{\ell^p_B} \, \leq \,
\left\| F \right\|_{L^p_B} \, \leq \,
A \left\| f \right\|_{\ell^p_B} \, .
\end{equation}
Here $A$ is a constant that depends only on $d${\rm ,} but not $p$ or the space $B${\rm .}
\endproclaim

Ideas of this kind go back to Plancherel and P\'{o}lya [PP].
In that work (when e.g.\ $d = 1$), the function $\frac{\sin \pi x}{\pi x}$ was used in effect in place of 
$\left(
\frac{\sin \pi x}{ \pi x} \right)^2$.  The resulting version of (2.3) is then more
delicate and holds only in the range $1 < p < \infty$, since it involves the Hilbert transform; 
it also does not cover the case of Banach space-valued functions.

To prove the lemma we observe two easily established estimates,
$$
\displaystyle{\int_{\Bbb{R}^d}} \, 
\Psi ( x ) d x \, \leq \,
A_1 \ \ {\rm and} \ \ \sup_{x} \:
\displaystyle{\sum_{n \in {\Bbb{Z}^d}}} \:
\Psi ( x - n ) \, \leq \, A_1\, .
$$
Then for any $p < \infty$, by H\"{o}lder' inequality
$$
\left| f_{\rm ext} ( x ) \right|^p \, \leq \,
\left(
\displaystyle{\sum_{n}} \, \left|
f ( n ) \right|^p \, \Psi ( x - n ) \right) \;
\left( \displaystyle{\sum_{n}} \, \Psi ( x - n ) \right)^{p - 1}
$$
and integration in $x$ then gives
$$
\parallel f_{\rm ext} \parallel^p_{L^p_B} \, \leq \,
\parallel f \parallel^p_{\ell^p_B} \, A^p_1 \, .
$$
The proof of the corresponding result for $p = \infty$ is similar but 
simpler.

To prove the converse inequality, choose $\widehat{\Phi}$ 
to be a $C^\infty$ function
with compact support so that $\widehat{\Phi} ( \xi ) = 1$, 
when $\xi \in 2 Q$.
Since 
$\left[ \left( \frac{\sin \pi x_1}{\pi x_1} \right)^2 \right]^\wedge \, = \,
( 1 - | \xi_1 | )_+$, it follows that 
$$
\widehat{\Psi}  \, \cdot \,
\widehat{\Phi} \, = \,
\widehat{\Psi} \, , \ {\rm and \ hence} \ \Psi \star \Phi \, = \, \Psi \, .
$$
As a result $f_{\rm ext} \star \Phi = f_{\rm ext}$ and since $f_{\rm ext} ( n ) = f ( n )$, we have
$$
f ( n ) \, = \,
\displaystyle{\int_{\Bbb{R}^d}} \,
f_{\rm ext} ( y ) \, \Phi ( n - y ) \, dy \, .
$$
Thus as before, 
$$
\big| f ( n ) \big|^p \, \leq \,
\left(
\displaystyle{\int_{\Bbb{R}^d}} \,
\big| f_{\rm ext} ( y ) \big|^p \,
\big| \Phi ( n - y ) \big| dy \right) \,
\left(
\displaystyle{\int_{\Bbb{R}^d}} \big|
\Phi ( n - y ) \big| \, dy \right)^{p - 1} \, ,
$$ 
and
$$
\left\| f \right\|^p_{\ell^p_B} \, \leq \,
\left\| f_{\rm ext} \right\|_{L^p_B} \, A^p_2 \, ,
$$
if
$$
\displaystyle{\int_{\Bbb{R}^d}} \, \left|
\Phi ( y ) \right| dy \, \leq \, A_2 \, , \ \ {\rm and} \ \
\sup_{y} \;
\displaystyle{\sum_{n \in {\Bbb{Z}^d}}} \;
\left| \Phi ( n - y ) \right| \, \leq \, A_2 \, . 
$$
The argument also gives the case $p = \infty$.
The lemma, inequality (2.3), is thus established with $A = \max ( A_1 , A_2)$.
 
To prove the proposition, we consider the cube $3 Q$ which can be covered
by $3^d$ disjoint translates of $Q$.  In fact, it is easily verified that  
$3Q \, = \, {\displaystyle\mathop{U}_{\ell \in {\Bbb{Z}^d , \ell = ( \ell_1 , \ldots \ell_d )}
\atop{\sup\limits_{j}} | \ell_j | \leq 1}} \, ( Q \, + \ell )$.
\noindent
Now let $m ( \xi )$ be continued periodically to $3 Q$, i.e. define $ \widetilde{m} ( \xi )$ by
$\widetilde{m} ( \xi ) = 
\displaystyle{\sum_{{\sup_j} | \ell_j | \leq 1}}$
$m ( \xi + \ell )$.  We let $\widetilde{T}$ 
denote Fourier multiplier operator, whose multiplier is $\widetilde{m} ( \xi )$.  
Then clearly
\begin{equation}
\parallel \widetilde{T} \parallel_{L^p_{B_1} 
\rightarrow L^p_{B_2}} \, \leq \,
3^d \, \parallel T \parallel_{L^p_{B_1} \rightarrow L^p_{B_2}}\, .
\end{equation}

On the other hand, we claim that
\begin{equation}
\widetilde{T} (f_{\rm ext}) \, = \,
\left( T_{\rm dis} ( f ) \right)_{\rm ext}\, .
\end{equation}
To verify (2.5) it suffices to do it for $f = \delta_m$, for every fixed $m$,
where 
$$ 
\delta_m ( n ) = \left\{
\begin{array}{lcl}
1 & {\rm if} & n = m \\
\\
0 & {\rm if} & n \ne m \, .
\end{array}
\right.
$$

We will check this by taking the Fourier transform of both
sides of (2.5).
Indeed 
$$ T_{\rm dis} ( f ) ( n ) \, = \, K ( n - m ) $$
and
$$
( T_{\rm dis} ( f ) )_{\rm ext} \, = \, 
\displaystyle{\sum_{n \in {\Bbb{Z}}^d}} \,
K ( n - m ) \: \Psi ( x - n )\, .
$$
Hence, 
$$
\begin{array}{ l c l}
\left(
T_{\rm dis} ( f )_{\rm ext} \right)^\wedge & = &
\displaystyle{\sum_{n}} \:
K ( n - m ) \;
\widehat{\Psi} \, ( \xi ) \, e^{- 2 \pi i n \xi} \\
\\
& = & \left(
\sum \: K ( n - m ) \, e^{- 2 \pi i n \xi} \right) \;
\widehat{\Psi} (  \xi ) \\
\\
& = & \left(
\displaystyle{\sum_{n}} \: K ( n ) \, e^{- 2 \pi i n \xi} \right)
\,
e^{- 2 \pi i m \xi} \;
\widehat{\Psi} ( \xi )\\ \\
& =& \, m_{\rm per} ( \xi) \,
\widehat{\Psi} ( \xi ) \,
e^{- 2 \pi i  m \xi} \hskip1.75in \hbox{(by (2.0)).}
\end{array}
$$
On the other hand
$$
\widetilde{T} ( f_{\rm ext} )^\wedge \, = \,
\tilde{m} ( \xi) \, 
( f_{\rm ext})^\wedge \, = \,
\widetilde{m} ( \xi) \, (\Psi ( x - m ))^\wedge \, = \,
\widetilde{m} ( \xi ) \widehat{\Psi} ( \xi) \, e^{- 2 \pi i m \xi} . $$ 
Now we have the desired identity,
since $\tilde{m} ( \xi) = m_{\rm per} ( \xi)$ on the support of 
$\widehat{\Psi}$ (note that $2 Q \subset 3 Q)$.

 \def\ds#1{{\displaystyle#1}}
Once (2.5) is established we have
$$
\begin{tabular}{rll}
 $\ds{\bigg\| T_{\rm dis} ( f ) \bigg\|_{\ell^p_{B_2}}}$ &$\leq 
\ds{A \bigg\| \left(
T_{\rm dis} ( f ) \right)_{\rm ext} \bigg\|_{L^p_{B_2}} }$&\qquad
 \hbox{(by the  lemma)}\\
&$=  \ds{A \bigg\| \widetilde{T} ( f_{ \rm ext} ) \bigg\|_{L^p_{B_2}} \, 
\leq \, 3^d \, A \bigg\| T ( f_{\rm ext} ) \bigg\|_{L^p_{B_2}}}$\\
&$\leq \ds{3^d \, A \bigg\| T \bigg\|_{L^p_{B_1} \rightarrow L^p_{B_2}} \,
\bigg\| f_{ \rm ext} \bigg\|_{L^p_{B_1}}}$&\qquad  \hbox{(by 
(2.4)))}\\
&$\leq \ds{3^d \, A^2 \bigg\| T \bigg\|_{L^p_{B_1} 
\rightarrow L^p_{B_2}} \,
\bigg\| f \bigg\|_{\ell^p_{B_1}}}$&\qquad  \hbox{(by  the   lemma).}\end{tabular}
$$
Thus $\bigg\| T_{\rm dis} \bigg\|_{\ell^p_{B_1} 
\rightarrow \ell^p_{B_2}} \, \leq \,
3^d \, A^2 \, \bigg\| T \bigg\|_{L^p_{B_1} \rightarrow L^p_{B_2}}$, and
the proposition is proved with $C = 3^d \, A^2$. 
\medbreak

We now fix an integer $q \geq 1$.  We shall also make
the stronger assumption that $m ( \xi)$ is supported in
$Q / q$, and consider $m^q_{\rm per}$ defined by
\begin{equation}
m^q_{\rm per} ( \xi ) \, = \,
\displaystyle{\sum_{\ell \in {\Bbb{Z}}^d}} \;
m ( \xi - \ell / q ) \, .
\end{equation} 

Notice that $m^q_{\rm per}$ is periodic with respect to elements in
$\left( 1 / q \right) {\Bbb{Z}^d}$, and hence, in particular, periodic
with respect to ${\Bbb{Z}^d}$.

We consider the operator $T^q_{\rm dis}$, a convolution operator
on ${\Bbb{Z}^d}$, having $m^q_{\rm per}$ as its Fourier multiplier; i.e.,
$$
\displaystyle{\sum_{m \in {\Bbb{Z}^d}}} \;
T^q_{\rm dis} \,
( f ) ( n ) \,
e^{- 2 \pi i n \xi} \, = \,
m^q_{\rm per} ( \xi ) \;
\displaystyle{\sum_{n \in {\Bbb{Z}^d}}} \,
f ( n ) \,
e^{- 2 \pi i n \cdot \xi} \, ,
$$ 
for suitable $f$.

\specialnumber{2.1}\proclaim{{C}orollary}
\begin{equation}
\bigg\| T^q_{\rm dis} \bigg\|_{\ell^p_{B_1} \rightarrow \ell^p_{B_2}} 
\, \leq \, C \bigg\| T \bigg\|_{L^p_{B_1} \rightarrow L^p_{B_2}}\, .
\end{equation}
Again the bound $C$ does not depend on $p${\rm ,} $B_1$ and $B_2${\rm ;} it is also
independent of~$q${\rm . }
\endproclaim

{\it Proof of the corollary}.
Let $T^q$ be the operator on
$L^p_{B_1} ( {\Bbb{R}^d} )$ to
$L^p_{B_2} ( {\Bbb{R}^d} )$ whose
multiplier is $m \left( \xi / q \right)$.  Notice that
$m ( \xi / q )$ is supported for $\xi \in Q$.

Now a simple scaling argument shows
$$
\bigg\|T^q \bigg\|_{L^p_{B_1} \rightarrow L^p_{B_2}} \, = \,
\bigg\| T \bigg\|_{L^p_{B_1} \rightarrow L^p_{B_2}} \, ,
$$
and so if $(T^q)_{\rm dis}$ is the discrete analogue in the
sense of Proposition 1, 
\begin{equation}
\bigg\| (T^q)_{\rm dis} \bigg\|_{\ell^p_{B_1} \rightarrow \ell^p_{B_2}} 
\, \leq \,
C \bigg\| T \bigg \|_{L^p_{B_1} \rightarrow L^p_{B_2}}\, .
\end{equation}

However, we must emphasize that $(T^q)_{\rm dis} \ne T^q_{\rm dis}$. 
In fact, the convolution kernel of $( T^q)_{\rm dis}$, which comes
from the multiplier $m ( \xi/q )$,  is $K^q ( n ) \, = \,
q^d K ( q n )$, $n \in {\Bbb{Z}^d}$.  
 
Next we observe the
convolution kernel, $K^\# ( n )$, of $T^q_{\rm dis}$ is given by
$$
K^{\#} ( m ) \, = \,
\left\{
\begin{array}{ll}
q^d K ( m ) , & {\rm if} \ m \in q {\Bbb{Z}^d} \\
\\
= \, 0 & {\rm if } \ m \in {\Bbb{Z}^d} \, , {\rm but} \ m \notin q {\Bbb{Z}^d} \, ,
\end{array}
\right.
$$
because
\begin{eqnarray*}
&&\displaystyle{\int_{Q}} \,
\left( 
\displaystyle{\sum_{\ell \in {\Bbb{Z}^d}}} \;
m ( \xi - \ell/q ) \right) \,
e^{2 \pi i \xi m} \, d \xi \\
&&\qquad = \; \displaystyle{\int_{Q}} \,
\left( 
\displaystyle{\sum_{\ell^\prime \in {\Bbb{Z}^d}}} \;
m ( \xi - \ell^\prime ) 
\right) \;
e^{2 \pi i \xi m} \, d \xi \, \times 
\left\{
\begin{array}{lcl}
q^d & {\rm if} & m \in q {\Bbb{Z}^d} \\
\\
0 & {\rm if} & m \notin q {\Bbb{Z}^d}\, .
\end{array}
\right.
\end{eqnarray*}
Now finally let $T_{\#}$ denote the operator mapping functions
of $q {\Bbb{Z}^d}$ to itself, given by the kernel $K^{\#}$, i.e.
$$
T_{\#} ( f^\prime) ( n q ) \, = \,
\displaystyle{\sum_{m \in {\Bbb{Z}^d}}} \,
f^\prime ( ( n - m ) q ) \,
K^{\#} ( m q ) \, .
$$

Then clearly
$$
\bigg\| T_{\#} \bigg\|_{{\ell^p_{B_1}} 
( q {\Bbb{Z}^d)} \, \rightarrow \, \ell^p_{B_2} ( q {\Bbb{Z}^d)}} \, = \,
\bigg\| ( T ^q )_{\rm dis} \bigg\|_{\ell^p_{B_1} ({\Bbb{Z}^d}) \, \rightarrow \,
\ell^p_{B_2} ( {\Bbb{Z}^d})} $$
which is an immediate consequence of the isomorphism 
${\Bbb{Z}^d} \, \leftrightarrow \, q {\Bbb{Z}^d}$, given by 
$n \leftrightarrow qn$, $( n \in {\Bbb{Z}^d})$.

Finally note that $T^q_{\rm dis}$ can be written as 
$T_{\#} \otimes I$, if we write $\ell^p_B ( {\Bbb{Z}^d} )$ \break as 
$\ell^p_B ( q {\Bbb{Z}^d} ) \otimes 
\ell^p ( {\Bbb{Z}^d}\diagup q {\Bbb{Z}^d} )$, 
with $T_{\#}$ acting on the first factor, and the identity acting on the 
second factor.

As a result
$$
\bigg\| T^q_{\rm dis} \bigg\|_{\ell^p_{B_1} \, \rightarrow \,
\ell^p_{B_2}} \, \leq \,
\bigg\| T_{\#} \bigg\|_{\ell^p_{B_1} ( q {\Bbb{Z}^d}) \, \rightarrow \,
\ell^p_{B_2} ( q {\Bbb{Z}^d}) } \, . 
$$ 
Combining this with (2.8) proves   Corollary  2.1.
\medbreak
 
We next consider a version of a convolution operator, whose multiplier is
somewhat akin to (2.4).  Here we shall consider 
\begin{equation}
m ( \xi ) \, = \,
\displaystyle{\sum_{\ell \in {\Bbb{Z}^d} }} \, \gamma_\ell \,
\Phi ( \xi - \ell/q )
\end{equation}
under the following assumptions:

\begin{namelist}{xxxxx}
\item [(a)] $\Phi$ is a $C^\infty$ function supported on $Q/q$.  As a function on $Q$ it has the Fourier expansion
$$
\Phi ( \xi) \, = \, 
\displaystyle{\sum_{m \in {\Bbb{Z}^d}}} \;
\varphi_m \, e^{- 2 \pi i m \xi}
$$
with
$$
\displaystyle{\sum_{m \in {\Bbb{Z}^d}}} \, | \varphi_m | \, \leq \,
A \, .
$$

\item[(b)] 
$\{ \gamma_\ell \}$ is a  $q {\Bbb{Z}^d}$ periodic sequence; i.e.,
$\gamma_\ell = \gamma_{\ell^\prime}$ if $\ell - \ell^\prime \in q {\Bbb{Z}^d}$.
\end{namelist}

Now let $\{ \hat{\gamma}_s \}$ be the Fourier transform of 
$\{ \gamma_\ell \}$; i.e., $\hat{\gamma}_s \, = \, 
\displaystyle{\sum_{\ell \in {\Bbb{Z}^d} \diagup q {\Bbb{Z}^d}}} \,
e^{ 2 \pi i s \ell/q} \, \gamma_\ell \, .$ 

 \vglue6pt
We shall also restrict our attention to scalar-related functions on
$\ell^p ( {\Bbb{Z}^d})$, as opposed to the Banach-space case treated in the
previous proposition, because of the specific use of Plancherel's identity.  
Our result is as follows.

\proclaim{Proposition}
Let $T$ be the operator on functions on ${\Bbb{Z}^d}$ whose 
Fourier multiplier is given by {\rm (2.8),} satisfying the conditions above{\rm .} 
Then with $1 \leq p \leq 2${\rm ,}
\begin{equation}
\bigg\| T \bigg\|_{\ell^p ( {\Bbb{Z}^d}) \, \rightarrow \,
\ell^p ( {\Bbb{Z}^d})} \, \leq \,
A \left( \sup_\ell \, | \gamma_\ell \,| \right)^{2 - 2 /p} \;
\left( \sup_s \, | \hat{\gamma}_s \, | \right)^{2 / p { -1}} \, .
\end{equation} 
\endproclaim

For the case $p = 2$ we use Parseval-Plancherel's theorem for 
${\Bbb{Z}^d}$, 
together with the disjointedness of the supports of the $\Phi ( \xi - \ell/q )$
and the fact that 
$\displaystyle{\sup_\xi} | \Phi ( \xi ) | \, \leq \, A$.  
This implies that
$ |m ( \xi ) | \, \leq A \, \displaystyle{\sup_{\ell}} | \gamma_\ell |$, 
yielding the case $p = 2$.

For the case $p = 1$, we calculate the $\ell^1 ( {\Bbb{Z}^d})$ norm
of the kernel $K ( n )$, corresponding to the multiplier $m ( \xi )$.  It
is given by 
\begin{eqnarray*}
K ( n )& = &
\displaystyle{\int_{Q}} \,
\left(
\sum \gamma_\ell \Phi ( \xi - \ell/q ) \right) \,
e^{2 \pi i n \xi} \, d \xi \\[5pt]
& = &
\varphi ( n ) 
\left( \displaystyle{\sum_{\ell \in {\Bbb{Z}^d} \diagup q {\Bbb{Z}^d}}} \;
\gamma_\ell \, e^{ 2 \pi i n \ell/q} \right) \, = \,
\varphi ( n ) \hat{\gamma}_n \, .
\end{eqnarray*}
Hence by property (a), 
$\displaystyle{\sum_{n \in {\Bbb{Z}^d}}} | K ( n )  | \, \leq \, A \,
\sup_n \, | \hat{\gamma}_n \,|$, and as a result the case 
$p = 1$ of (2.10) is proved.  The general result for $1 \leq p \leq 2$ then
follows by Riesz' convexity theorem.  

 \section{The main term}

The averages we are interested in, 
$$
A_\lambda ( f ) ( n ) \, = \,
\frac{1}{N ( \lambda )} \,
\displaystyle{\sum_{| m | = \lambda} } \, f ( n - m ) \, ,
$$
will be replaced by the equivalent averages
$$
\frac{1}{\lambda^{d - 2}} \,
\displaystyle{\sum_{| m | = \lambda}} \,
f ( n - m ) \, ,
$$
when $d \geq 5$.  This equivalence comes about because, as we have
pointed out with 
$N ( \lambda )$ = number of 
$n \in {\Bbb{Z}}^d$, so that 
$|n| = \lambda$, we have 
$N ( \lambda ) \approx \lambda^{d - 2}$, 
whenever $\lambda^2$ is an integer, and 
$d \geq 5$.  In order not to introduce new notation, we shall
designate these averages also by $A_\lambda$ and now write
\setcounter{equation}{0}
\begin{equation}
A_\lambda ( f ) ( n ) \, = \,
\frac{1}{\lambda^{d - 2}} \,
\displaystyle{\sum_{| m| = \lambda}} \,
f ( n - m ) \,, 
\end{equation}
and  in what follows we shall always  assume that $\lambda$ is restricted 
so that $\lambda^2$ is an integer.

Here we shall deal with the main term in the approximation of
$A_\lambda$.  It is a convolution operator $M_\lambda$
acting on functions on ${\Bbb{Z}}^d$, which can be written as
\begin{equation}
M_\lambda \, = \,
c_d \,
\displaystyle{\sum^{\infty}_{q = 1}} \,
\displaystyle{\sum_{1 \leq a \leq q \atop {(a , q ) = 1}}} \,
e^{- 2 \pi i \lambda^2 a / q} \, M^{a/q}_\lambda \, ,
\end{equation} 
where the sum is taken over all reduced fractions $a/q$, with
$0 < a/q \leq 1$.  Hence $c_d$ is the constant =
$\frac{\pi^{d/2}}{\Gamma (d/2)}$.  Also each $M^{a/q}_\lambda$ 
is the convolution operator
whose multiplier is 
\begin{equation}
\displaystyle{\sum_{\ell \in {\Bbb{Z}^d}}} \,
G ( a / q , \ell ) \,
\Psi_q \, (\xi - \ell/q) \,
 d \widehat{\sigma}_\lambda \, ( \xi - \ell/q ) \, .
\end{equation}

In the above, $\Psi_q ( \xi ) = \Psi (q \xi )$, where $\Psi$ is a 
$C^\infty$ cut-off function supported in the cube
$Q/2$, with $\Psi ( \xi ) = 1$, for $\xi \in Q/4$.  Also 
$G ( a / q , \ell )$ is the normalized Gauss sum
 \begin{equation}
G ( a / q , \ell ) \, = \,
q^{-d} \,
\displaystyle{\sum_{n \in {\Bbb{Z}^d}/q {\Bbb{Z}^d} }} \,
e^{2 \pi i ( |n|^2 a / q + n \cdot \ell/q )} \, ,
\end{equation} 
and $d \hat{\sigma}_\lambda ( \xi ) $ is the Fourier transform of the
normalized invariant measure $d \sigma_\lambda$ 
supported on the
sphere $|x| = \lambda$.  Note that (3.3) is periodic on 
$\xi$ with periods in ${\Bbb{Z}^d}$, since 
$G ( a / q , \ell ) = G(a/q , \ell^\prime )$, if 
$\ell \equiv \ell^\prime {\rm mod} \, q {\Bbb{Z}^d}$; also, for each $\xi$
only one term in (3.3) is nonzero.

We define the corresponding maximal operators, 
$$M_\star ( f ) ( n ) \, = \,
\displaystyle{\sup_{0 < \lambda < \infty}} \,
\left|
M_\lambda ( f ) ( n ) \right|,$$
and 
$$M^{a / q}_\star ( f ) ( n ) \, = \,
\displaystyle{\sup_{0 < \lambda < \infty}} \, 
\left| M_\lambda^{a / q} ( f ) ( n ) \right|.$$

The basic estimates for these are as follows:

\proclaim{Proposition} 
{\rm (a)} 
$\parallel M^{a/q}_\star 
\parallel_{\ell^p \rightarrow \ell^p} \, = \, 
O \left( q^{-d ( 1 - 1/p )} \right)$  
{if} $d \geq 3${\rm ,} \ {and} \ $\frac{d}{d-1} < p \leq 2,$
\medbreak
{\rm (b)} 
$\parallel M_\star
\parallel_{\ell^p \rightarrow \ell^p} \, \leq \, A $ 
 if $d \geq 5${\rm ,} and   $\frac{d}{d - 2} \, < \, p \, \leq \, 2 \, .$
\endproclaim
 
To prove part (a) we write 
$\Psi = \Psi \cdot \Psi^\prime$, where 
$\Psi^\prime$ is another $C^\infty$ function, supported in $Q$, with
$\Psi^\prime ( \xi) = 1$
for $\xi \in Q/2$.  Then the operator corresponding to the
multiplier (3.3) can be written as a product of two
operators, with multipliers respectively:
$$
\displaystyle{\sum_{\ell \in {\Bbb{Z}^d}}} \;
G ( a / q , \ell ) \,
\Psi^\prime_{q} \, ( \xi - \ell/q )
$$
and
$$
\displaystyle{\sum_{\ell \in {\Bbb{Z}^d}}} \;
\Psi_q ( \xi - \ell/q ) \,
d \hat{\sigma}_\lambda \, ( \xi - \ell/q ) \, ,
$$
where $\Psi^\prime_q ( \xi ) \, = \, \Psi^\prime ( q \xi )$, if we recall
that for each $\xi$ only one term in each of the above sums is 
nonvanishing.

To the first multiplier we apply Proposition 2.2 (in $\S$2) with 
$\gamma_\ell = G ( a / q , \ell )$, and 
$\Phi ( \xi ) = \Psi^\prime_q ( \xi )$.  Note that
$\Phi ( \xi )$ is supported in $Q/q$, and moreover
$$
\varphi ( n ) \, = \,
\displaystyle{\int_{{\Bbb{R}^d}}} \,
\Phi ( \xi ) \, e^{- 2 \pi i n \xi} \, d \xi \, = \,
q^{-1} \, \widetilde{\Psi} ( q^{-1} n ) \, ,
$$
where $\widetilde{\Psi}$ is the Fourier transform of $\Psi^\prime$. 
Now $| \widetilde{\Psi} ( x ) \, | \leq \, 
A_N ( 1 + | x |)^{-N}$
for all $N \geq 0$, so that 
$\displaystyle{\sum_{n \in {\Bbb{Z}^d}}} \, | \varphi ( n ) | \leq A$.

Next, there is the estimate 
$| G ( a/q , \ell ) | \, = \,
O ( q^{-d/2} )$;  this is well-known, but in any case it follows  from the
standard one-dimensional case merely by observation that
$G ( a / q , \ell )$ is a $d$-fold product of these one-dimensional 
sums.\footnote{For the one-dimensional estimates, see [W].}  
Moreover, if 
$\hat{\gamma}_s  = 
\displaystyle{\sum_{\ell \in {\Bbb{Z}^d}/q {\Bbb{Z}^d}}} \,
e^{2 \pi i s \cdot \ell/q} \, G(a/q , \ell )$, then,
$$
\hat{\gamma}_s \, = \,
\frac{1}{q^d} \,
\displaystyle{\sum_{n \in {\Bbb{Z}^d/_{q {\Bbb{Z}^d}}}}} \;
\displaystyle{\sum_{\ell \in {\Bbb{Z}^d/q} \, {\Bbb{Z}^d}}} \;
e^{2 \pi i s \cdot \ell/q} \,
e^{2 \pi i (a/q ) | n |^2} \,
e^{2 \pi i n \cdot \ell/q} \, = \,
e^{2 \pi i (a/q) | s|^2} \, .
$$
Hence by Proposition 2.2, the norm of the corresponding operator
(acting on $\ell^p$ to itself, $1 \leq p \leq 2$, with
scalar-valued functions) is
$O ( q^{- (d/2) ( 2- 2/p)} ) = O ( q^{-d (  1 \, - 1/p)} )$.

Next, the multiplier
\begin{equation}
\displaystyle{\sum_{\ell \in {\Bbb{Z}^d}}} \,
\Psi_q \, ( \xi - \ell/q ) \, d \hat{\sigma}_\lambda ( \xi - \ell/q )
\end{equation}
 corresponds to a convolution operator from
$\ell^p ( {\Bbb{Z}^d} )$  (scalar-valued), to 
$\ell^p_B ( {\Bbb{Z}^d} )$, where $B$ is the $\ell^\infty$ 
space of functions of $\lambda > 0$, for which $\lambda^2$ is an
integer, and $0 < \lambda^2 \leq N$.  Notice that
$\Phi_q ( \xi ) = \Phi ( q \xi )$ is a bounded multiplier of 
$L^p ( R^d )$ to itself (with norm 
independent of $q$).  Observe also that $d \hat{\sigma}_\lambda ( \xi )$
is a bounded multiplier from
$L^p ( R^d )$ to $L^p_B ( {\Bbb{R}^d} )$, for
$p > \frac{d}{d-1}$, which is a consequence of the spherical maximal 
theorem in ${\Bbb{R}^d}$.  Finally, note that
$m ( \xi ) = \Phi_q ( \xi ) \,
d \hat{\sigma}_\lambda ( \xi )$ is supported in $Q/q$.  Thus,
applying the corollary to Proposition 2.1, we see that (3.5) is a 
bounded multiplier from $\ell^p ( {\Bbb{Z}^d})$ to 
$\ell^p_B ( {\Bbb{Z}^d})$, with norm independent of $N$  (and $q$).   
Letting $N \rightarrow \infty$, and combining this with the estimate for
the first multiplier, we have established 
conclusion (a) of Proposition 3.1.  The second conclusion
follows from this because 
$$
M_\star \leq c_d \;
\displaystyle{\sum_{1 \leq q < \infty}} \;
\displaystyle{\sum_{(a , q ) = 1 \atop {1 \leq a \leq q}}} \;
M^{a/q}_\star \, , $$
so that
$$
\bigg\| M_\star \bigg\|_{\ell^p \rightarrow \ell^p} \, \leq \,
A \,
\displaystyle{\sum_{1 \leq q < \infty}} \,
q \cdot q^{-d ( 1 - 1/p)} \, < \, \infty \, ,
$$
 if $1 - d ( 1 - 1/p ) < - 1$,
  i.e.\  when $p > \frac{d}{ d - 2  }$.
 
\section{Approximations}

We now state the assertions which guarantee that $M_\lambda$
provides an adequate approximation to our operator $A_\lambda$.
There are two facts; the first is a purely $\ell^2$ statement.

\proclaim{Proposition}
There is a bound $A${\rm ,} so that for any $\Lambda > 0${\rm ,} 
\setcounter{equation}{0}
\begin{equation}
\parallel 
\displaystyle{\sup_{\Lambda \leq \lambda \leq 2 \Lambda}} \;
| A_\lambda ( f ) \, - 
M_\lambda ( f ) | \,
\parallel_{\ell^2} \, \leq \,
A \, \Lambda^{2 - d /2} \parallel \, f \parallel_{\ell^2}\, ,  \ \ 
\hbox{if }  d  \geq  5 \, .\hskip.45in
\end{equation} 
\endproclaim

The second is a partial result for $A_\lambda$ which was known
previously (see~[M]).

\proclaim{Proposition}
There is a bound $A${\rm ,} so that for any $\Lambda > 0$
\begin{equation}
\parallel \;
\displaystyle{\sup_{\Lambda \leq \lambda \leq 2 \Lambda}} \;
| A_\lambda ( f ) \, \parallel_{\ell^p} \, \leq \,
A \parallel \, f \, \parallel_{\ell^p} \, , \ \  
\hbox{if } \ d \geq 5 \, , \hspace{.8em} p > \frac{d}{d-2} \, .\hskip.5in
\end{equation}
\endproclaim

Recall that the $\lambda$ which appear in (4.1) and (4.2) are always restricted
to the fact that $\lambda^2$ is an integer.

We shall momentarily take these two propositions for granted and
see how they, together with Proposition (3.1), prove our main theorem.

Now (4.2) together with Proposition (3.1) yield
$$
\parallel \displaystyle{\sup_{\Lambda \leq \lambda \leq 2 \Lambda}} \;
(A_\lambda - M_\lambda ) f \, \parallel_{\ell^p} \, \leq \,
A \, \parallel \, f \, \parallel_{\ell^p} \hskip.5in 
\hbox{\rm for }   2   \geq  p   >  \frac{d}{ d - 2} \, .
$$
Interpolating this with (4.1) gives 
$$
\parallel 
\displaystyle{\sup_{\Lambda \leq \lambda \leq 2 \Lambda}} \;
\left| ( A_\lambda \, - M_\lambda ) ( f ) \right| \,
\parallel_{\ell^p} \, \leq \,
A \, \Lambda^{- \varepsilon ( p )} \,
\parallel \, f \, \parallel_{\ell^p}
$$
for some $\varepsilon ( p ) > 0$, if $\frac{d}{d - 2} < p \leq 2$.

Next, 
$$
\displaystyle{\sup_{1 \leq \lambda < \infty}} \, \left|
 ( A_\lambda \, - M_\lambda ) \,  ( f )  \right| \, \leq \,
\displaystyle{\sum^{\infty}_{k = 0}} \;
\displaystyle{\sup_{2^k \leq \lambda \leq 2^{k+1}}} \,
\left| ( A_\lambda \, - M_\lambda ) f \right|.
$$
Taking  the $\ell^p$ norm we get that 
$$
\parallel \,
\displaystyle{\sup_{1 \leq \lambda < \infty}} \, \left|
( A_\lambda \, - M_\lambda ) \,  f \right| \, 
\parallel_{\ell^p} \, \leq \,
A^\prime \, \parallel \, f \, \parallel_{\ell^p} \, ,
\ {\rm for} \ \frac{d}{d - 2} \, < \, p \, \leq \, 2 \, ,
$$
since $\displaystyle{\sum_{k}} \, 
2^{- \xi ( p ) k} < \infty$.  Thus, invoking
Proposition (3.1) again yields
$$ 
\parallel \, 
\displaystyle{\sup_{1 \leq \lambda < \infty}} | 
\, A_\lambda ( f ) | \parallel_{\ell^p} \leq A \parallel \, f \, 
\parallel_{\ell^p} \quad \hbox{when} \quad 
\frac{d}{d - 2} \, < \, p  \,\leq \, 2 \, .
$$ 
 
Since the corresponding estimate for $p = \infty$ is trivial, the full
range $\frac{d}{d - 2} < p \leq \infty$ then follows by
interpolation, proving the main theorem.

\section{The decomposition of $A_\lambda$}

\hspace{.25in} 
To prove the crucial approximation property (4.1) we shall decompose the
operator $A_\lambda$ into a sum, each of whose terms corresponds to
a fraction $a / q$, with $1 \leq q$, 
$1 \leq a \leq q$, and $(a , q ) = 1$.  It is here we use the
ideas of the ``circle method'' of Hardy, Littlewood, and Ramanujan.

Let us fix $\Lambda > 0$, and consider any $\lambda$ for which
$\Lambda \leq \lambda \leq 2 \Lambda$.  We shall write $a_\lambda ( \xi )$
for the multiplier corresponding to the operator 
$A_\lambda$ given by (3.1).  We claim that
$$
a_\lambda ( \xi ) \, = \,
\frac{e^{2 \pi \varepsilon \lambda^2}}{\lambda^{d-2}} \;
\displaystyle{\sum_{n \in {\Bbb{Z}^d}}} \;
e^{- 2 \pi \varepsilon | n |^2} \,
e^{2 \pi i n \cdot \xi} \;
\displaystyle{\int^{1}_{0}} \;
e^{2 \pi i ( |n|^2 \, - \lambda^2 ) t} \; dt\, .
$$
Here $\varepsilon$ is positive, but otherwise arbitrary;  we will fix it
later by setting $\varepsilon = 1/\Lambda^2$.  This identity is obvious
because 
$\displaystyle{\int^{1}_{0}} \, 
e^{2 \pi i ( | n |^2 - \lambda^2)t} \, dt = 1$ or 
$0$ according to whether $|n| = \lambda$ or not.

Now we introduce the $\Theta$ function
\setcounter{equation}{0}
\begin{equation}
{\cal{F}} ( z , \xi ) \, = \,
\displaystyle{\sum_{n \in {\Bbb{Z}^d}}} \,
e^{- 2 \pi | n |^2 z} \,
e^{2 \pi i n \xi} \, ,
\end{equation}
for $\Re ( z ) > 0$, and we make a Farey direction of level = $\Lambda$
of the interval $[ 0 , 1]$ of the $t$ integration.  That is, for
each $a / q$, $(a , q ) = 1$ with
$1 \leq a \leq q$, and $q \leq \Lambda$, we associate the
interval 
$\bar{I} ( a / q ) = \left\{ t : \, - \frac{\beta}{q \Lambda} \, 
\leq \, t - a/q \, \leq \, \frac{\alpha}{q \Lambda} \right\}$, where 
$\alpha = \alpha ( a/q , \Lambda ) \approx 1$, and 
$\beta = \beta ( a/q , \Lambda ) \approx 1$, with $\alpha$ and
$\beta$ chosen appropriately.  We denote by $I ( a / q )$ the
corresponding intervals translated to the origin,
$I ( a / q ) = \left\{ \tau: \, 
- \frac{\beta}{q \Lambda} \, \leq \, \tau \, \leq \,
\frac{\alpha}{q \Lambda} \right\}$.  Inserting this in 
the above formula for $a_\lambda ( \xi)$ and using identity (5.1) we get 
$$
a_\lambda ( \xi) \, = \,
\displaystyle{\sum_{1 \leq q \leq \Lambda}} \;
\displaystyle{\sum_{1 \leq a \leq q \atop {(a , q ) = 1}}} \;
a_\lambda^{a / q } ( \xi ) 
\, ,
$$
where
\begin{equation}
a_\lambda^{a / q} ( \xi) \, = \,
\frac{e^{2 \pi \varepsilon \lambda^2}}{\lambda^{d - 2}} \,
e^{- 2 \pi i \lambda^2 a /q} \;
\displaystyle{\int_{I ( a / q )}} \;
e^{- 2 \pi i \lambda^2 \tau} \;
{\cal{F}} ( \varepsilon - i \tau - i \, a / q , \xi ) \, d \tau \, .\quad
\end{equation}

Next we use the fundamental identity for the $\Theta$ function (5.1).
It states that for $\Re ( z ) > 0$,
\begin{equation}
{\cal{F}} \left( z - i \, \frac{a}{q} , \xi \right) \, = \,
\frac{1}{(2z)^{d / 2}} \;
\displaystyle{\sum_{\ell \in {\Bbb{Z}^d}}} \;
G \left( \frac{a}{ q }, \ell \right) \: \exp 
\left( \frac{- \pi | \xi - \frac{\ell}{q} |^2}{2z} \right) \, .\hskip.5in
\end{equation}
Here $G ( a / q , \ell )$ is the normalized Gauss sum (3.4).  The above
is the $d$-dimensional version of a familiar identity.  (For $d = 1$ 
see, e.g.,  [SW$_1$, (3.4)]; also [W].)  
The general case $d \geq 1$ can be proved
the same way invoking the Poisson summation formula; alternatively one 
can observe that (5.3) is merely the $d$-fold product of the 
corresponding $1$-dimensional identities for each variable
$\xi_1 , \xi_2 , \ldots \xi_d$, separately.

From (5.3) and (5.2) it follows that
\begin{equation}
a^{a/q}_\lambda ( \xi ) \, = \,
e^{- 2 \pi i \lambda^2 a / q} \;
\displaystyle{\sum_{\ell \in {\Bbb{Z}^d}}} \;
G ( a / q , \ell ) \,
J_\lambda \, (a/q , \xi \, - \ell/q ) \, , 
\end{equation}
where
\begin{equation}
J_\lambda ( a / q , \xi ) \, = \,
\frac{e^{2 \pi \varepsilon \lambda^2}}{\lambda^{d - 2}} \;
\displaystyle{\int_{I ( a / q )}} \;
e^{- 2 \pi i \lambda^2 \tau} \,
( 2 ( \varepsilon \, - i \tau ))^{-d/2} \;
e^{ \frac{ - \pi | \xi |^2}{2 ( \varepsilon \, - i \tau )}} \,
d \tau \, .\hskip.4in
\end{equation}

\section{Approximations, continued}

We shall approximate the multipliers $a_\lambda^{a / q} ( \xi )$ above
by multipliers $b^{a/q}_{\lambda} ( \xi )$ where the cut-off factors
$\Phi_q ( \xi \, - \ell/q )$ have been inserted in (5.4). That is, we define 
\setcounter{equation}{0}
\begin{equation}
b^{a/q}_\lambda ( \xi ) \, = \,
e^{- 2 \pi i \lambda^2 a / q} \;
\displaystyle{\sum_{\ell \in {\Bbb{Z}^d}}} \;
G ( a / q , \ell ) \,
\Phi_q ( \xi \, - \ell/q ) \;
J_\lambda ( \xi \, - \ell/q ) \, .\hskip.4in
\end{equation}
Here $\Phi_q ( \xi ) = \Phi ( q \xi )$.

Next we approximate $b^{a / q}_\lambda ( \xi )$ by replacing the
integral (5.5) that appears in (6.1) by the corresponding integration
when taken over the whole real line.  So we set 
\begin{equation}
c^{a/q}_{\lambda} ( \xi ) \, = \,
e^{- 2 \pi i \lambda^2 a / q} \;
\displaystyle{\sum_{\ell \in {\Bbb{Z}^d}}} \;
G ( a / q , \ell ) \;
\Phi_q ( \xi \, - \ell/q ) \;
I_\lambda ( \xi \, - \ell/q ) \, ,\hskip.4in
\end{equation} 
with
\begin{equation}
I_\lambda ( \xi ) \, = \,
\frac{e^{2 \pi \varepsilon \lambda^2}}
{\lambda^{d - 2}} \;
\displaystyle{\int^{\infty}_{- \infty}} \;
e^{-2 \pi i \lambda^2 \tau} \;
( 2 ( \varepsilon \, - i \tau ) )^{-d/2} \;
e^{ \frac{- \pi |\xi |^2}{2 ( \varepsilon \, - i \tau )}} \, d \tau \, .
\end{equation}
 
We define the operators $A^{a/q}_{\lambda}$, $B^{a/q}_{\lambda}$, 
$C^{a/q}_{\lambda}$, as the convolution operators (acting on functions
of ${\Bbb{Z}^d}$), whose Fourier multipliers are
respectively, $a^{a/q}_\lambda ( \xi )$, $b^{a/q}_\lambda ( \xi )$, and
$c^{a / q}_\lambda ( \xi )$. 

\proclaim{Proposition}
\begin{eqnarray}
\qquad\displaystyle{\sum_{1 \leq q \leq \Lambda}} \;
\displaystyle{\sum_{1 \leq a \leq q \atop { (a , q ) = 1}}} \;
\parallel
\displaystyle{\sup_{\Lambda \leq \lambda \leq 2 \Lambda}} \:
\left| \,\left(
A^{a / q }_\lambda \, - B_\lambda^{a / q} \right) f \, \right| \,
\parallel_{\ell^2}&\hskip-7pt\leq\hskip-7pt&
A \, \Lambda^{2 - d/2} \,
\| \, f \, \|_{\ell^2}\, ,
\\
\displaystyle{\sum_{1 \leq q \leq \Lambda}} \;
\displaystyle{\sum_{1 \leq a \leq q \atop { (a , q ) = 1}}} \;
\parallel 
\displaystyle{\sup_{\Lambda \leq \lambda \leq 2 \Lambda}} \;
\left|
\left(
B^{a / q}_{\lambda} \, - \,
C^{a / q}_\lambda \right) \, f \right| \, \parallel_{\ell^2}&\hskip-7pt\leq\hskip-7pt& A \, \Lambda^{2 - d/2} \, \parallel \, f \, \parallel_{\ell^2} \, .
\end{eqnarray} 
\endproclaim
 
It is understood that in the above assertions our $\varepsilon$ is fixed
to be = $1 / \Lambda^2$.

To prove (6.4), let $F_\tau$ be the function on
${\Bbb{Z}^d}$ which is given in terms of its Fourier expansion
by $\hat{F}_\tau ( \xi ) = \mu ( \xi ) \, \hat{f}  ( \xi )$
where 
$$
\mu ( \xi ) = \displaystyle{\sum_{\ell \in {\Bbb{Z}}}} \;
( 1 \, - \Phi_q ( \xi \, - \ell/q ) ) \,
e^{- \pi | \xi \, - \ell/q |^2/2 ( \varepsilon \, - i \tau )} \, .
$$ 
Note that since each term in the sum is supported where
$| \xi \, - \ell/q |\, \geq \, c / q$,
$$\displaystyle{\sup_{\xi}} \, | \mu ( \xi ) | \, \leq \, A \, \exp 
\left( \frac{ - c \varepsilon}{q^2 ( \varepsilon^2 + \tau^2 )} \right),$$
for some $c > 0\, .$
Thus
$$
\parallel \, F_\tau \, \parallel_{\ell^2} \, \leq \,
A \, \exp 
\left(
\frac{- c \varepsilon}{q^2 ( \varepsilon^2 + \tau^2 )} \right) \parallel \, f \,
\parallel_{ \ell^2}\, .
$$ 

Now observe that
$$
{\displaystyle{\sup_{\Lambda \leq \lambda \leq 2 \Lambda}} \,
\left| \left( A^{a / q}_{\lambda} \, - 
B^{a / q}_\lambda \right) \, f \, \right|} \; 
\leq \;
A \, \Lambda^{-d + 2} \,
q^{- d / 2} \,
\displaystyle{\int_{I ( a / q )}} \;
( \varepsilon^2 + \tau^2 )^{-d/4} \, | F_\tau | d \tau \, ,
$$
because $G ( a / q , \ell ) \, = \,
O ( q^{- d / 2} )$.  

\vglue2pt
As a result,
\begin{eqnarray*}
\lefteqn {\displaystyle{\sup_{\Lambda \leq \lambda \leq 2 \Lambda}} \;
\left\| \, \left|\,
\left( A^{a / q}_{\lambda} \, - B^{a / q}_{\lambda} \right) \, f \, \right|
\right\|_{\ell^2} } \\[4pt]
& & \leq \; A \, \Lambda^{- d + 2} \,
q^{-d / 2} \;
\displaystyle{\int_{I ( a / q )}} \;
( \varepsilon^2 + \tau^2 )^{-d/4} \,
\exp 
\left( \frac{- c \varepsilon}{q^2 ( \varepsilon^2 + \tau^2 )} \right) \,
d \tau \, \cdot \, \parallel \, f \, \parallel_{\ell^2} \, .
\end{eqnarray*}
Now, because $e^u \leq c \, u^{-d/4}$, we get a contribution of
$\Lambda^{-d + 2} \, 
\varepsilon^{-d/4} \, | I ( a / q ) | \break \parallel \, f \, \parallel_{\ell^2}$.
Taking into account that $\varepsilon = 1/ \Lambda^2$, and
$\displaystyle{\sum_{1 \leq q \leq \Lambda}} \;
\displaystyle{\sum_{a ,q}} \; 
| I ( a / q ) | = 1$, we obtain (6.4).

The proof of (6.5) is similar.  Notice that we are now integrating
over $\tau$ in the complement of $I ( a / q )$, and thus
$| \tau | \geq c/q \Lambda$.  We are led in the same way to see that
\begin{eqnarray*}
\lefteqn{ \| \,
\displaystyle{\sup_{\Lambda \leq \lambda \leq 2 \Lambda}} \,
\left| \left( B^{a/q}_{\lambda} \, - \, 
C^{a/q}_\lambda \right) \, f \, \right| \,
\|_{\ell^2} } \\
\\
&& \leq \, A \, \Lambda^{-d + 2} \, q^{- d/2} \,
\displaystyle{\int_{| \tau | \geq c/q\Lambda}} \;
\tau^{- d/2} \, d \tau \, \cdot \, \parallel \, f \, 
\parallel_{\ell^2} \\
\\
&& \leq \,
A \, \Lambda^{-d + 2} \, q^{-d/2} \, ( q \Lambda )^{d/2 - 1} \,
\parallel \, f \parallel_{\ell^2} \, = \,
A \, \Lambda^{-d/2 + 1} \, q^{-1} \, \parallel \, f \, \parallel_{\ell^2}\, .
\end{eqnarray*} 
Now sum in $a$, then over $q$, $q \leq \Lambda$.  This gives
a contribution of 
$$
 O \left( 
\Lambda^{-d/2 + 1} \right) \; 
\left(
\displaystyle{\sum_{1 \leq q \leq \Lambda}} \hspace{.15in}
\displaystyle{\sum_{1 \leq a \leq q}} \, 1 
\right) \, 
q^{-1} \, 
\parallel \, f \parallel_{\ell^2} 
 \;  = \; O \, 
\left( 
\Lambda^{-d / 2 + 2} 
\right) \, 
\parallel \, f \, \parallel_{ \ell^2},
$$
which proves (6.5).
 
To complete the approximation process (the proof of (4.1)) we now
identify $I_\lambda ( \xi )$ given by (6.3).

\specialnumber{6.1}\proclaim{Lemma}
$$
I_\lambda ( \xi ) \, = \, c_d \, d \hat{\sigma}_\lambda \, ( \xi ) \, .
$$
\endproclaim

Taking this temporarily for granted we observe that as a result,
$c_d \, M^{a/q}_\lambda = C^{a/q}_\lambda$ (see (3.3) and (6.2),
(6.3)).  Hence for 
$\Lambda \leq \lambda \leq 2 \Lambda$,
\begin{eqnarray*}
\left|
\left(
A_\lambda \, - M_\lambda \right) \, f \, \right| &\leq &
\displaystyle{\sum_{1 \leq q \leq \Lambda}} \ \ 
\displaystyle{\sum_{(a , q ) = 1}} \;
\left|
\left(
A^{a/q}_\lambda \, - C^{a/q}_{\lambda} \right) \, f \, \right|
\\
&& + \ 
\displaystyle{\sum_{q > \Lambda}} \ \
\displaystyle{\sum_{(a , q ) = 1}}  \,
c_d \, \left| M^{a/q}_{\lambda} \, f \, \right| \, .
\end{eqnarray*}
However,
$$
\left| \,
A^{a / q}_\lambda \, - C^{a/q}_{\lambda} \, \right| \, \leq \, 
\left| \, A^{a/q}_\lambda \, - B^{a / q}_{\lambda} \, \right| \, + \,
\left| \, B^{a/q}_\lambda \, - C^{a/q}_\lambda \, \right| .
$$
Thus, invoking Proposition (6.1), and Proposition (3.1) for $p = 2$, we see that
\begin{eqnarray*}
 &&\hskip-30pt\parallel \,
\displaystyle{\sup_{\Lambda \leq \lambda \leq 2 \Lambda}} \, \left| \,
\left( 
A_\lambda \, - M_\lambda 
\right) \, ( f ) \,\right| \, \parallel_{\ell^2}   \\
& & =\  O \left(
\Lambda^{2 - d/2} \right) \, \left\| \, f \right\|_{\ell_2} \, + \,
\left(
\displaystyle{\sum_{q > \Lambda}} \ \ 
\displaystyle{\sum_{1 \leq a \leq q}} \
q^{-d/2} \right) \, \left\| \, f \right\|_{\ell^2} \\
\\
&& =\ O \left(
\Lambda^{2 - d /2} \right) \, \left\| \, f \, \right\|_{\ell^2} \, , \
{\rm if} \ \ d \, \geq \, 5 \, .
\end{eqnarray*}

Therefore, Proposition (4.1) is now proved, and with it the proof
of our main theorem is complete, save for   verification of the lemma
above. 

\medskip
\section{Proof of   Lemma 6.1}

The identity
\setcounter{equation}{0}
\begin{equation}
\frac{e^{2 \pi \varepsilon \lambda^2}}{\lambda^{d - 2}} \;
\displaystyle{\int^{\infty}_{- \infty}} \;
e^{- 2 \pi i \lambda^2 \tau} \;
( 2 ( \varepsilon \, - i \tau ) )^{-d/2} \;
e^{ \frac{- \pi | \xi |^2}{2 ( \varepsilon \,- i \tau )}} \, d \tau \, = \,
c_d \, d \hat{\sigma}_\lambda ( \xi )\hskip.5in
\end{equation}
is probably known, but we have not found it in the literature, and
so we will give a proof.

First we observe that the left-side of (7.1) is in fact
independent of $\varepsilon$, and so we may take $\varepsilon = 1 / \lambda^2$.
We see that this follows by changing the contour while integrating the function
$F(z) = (2z)^{-d/2} \, 
e^{2 \pi z} \, 
e^{- \pi | \xi |^{2}/2 z}$ along lines 
parallel to the $x$ axis in the upper half-plane.  Next, with
$\varepsilon = 1/\lambda^2$, and with the change of variables
$\lambda^2 \tau = t$, 
$$
I_\lambda ( \xi ) \, = \,
e^{2 \pi} \,
\displaystyle{\int_{- \infty}^{\infty}} \;
e^{- 2 \pi i t} \;
\frac{1}{(2 ( 1 \, - i t ))^{d/2}} \;
e^{ \frac{ - \pi \lambda^2 | \xi |^2}{2 ( 1 \, - i t )}} \, dt\, .
$$

We now insert an extra convergence factor $e^{- \pi \delta t^2}$ in the
integral defining $I_\lambda$ above.  Denoting the resulting integral
by $I^\delta_\lambda$ we have 
$I^\delta_\lambda \rightarrow I^\delta_\lambda$; moreover if 
$\varphi$ is any test function in the Schwartz space, then 
$$
\displaystyle{\int_{\Bbb{R}^d}} \,
\widehat{\varphi} ( \xi ) \,
I_\lambda ( \xi ) \, d \xi \, = \,
\lim_{\delta \rightarrow 0} \,
\displaystyle{\int_{\Bbb{R}^d}} \, \widehat{\varphi} ( \xi) \,
I^\delta_\lambda ( \xi) \, d \xi\, .
$$
Also,
\begin{equation}
\displaystyle{\int_{\Bbb{R}^d}} \,
\widehat{\varphi} ( \xi ) \,
I^\delta_{\lambda} ( \xi) \, d \xi \, = \,
\displaystyle{\int_{\Bbb{R}^d}} \,
\varphi ( x ) \, \hat{I}^\delta_\lambda ( x ) \, d x \, .
\end{equation}

Calculating the Fourier transform of the Gaussian
$e^{ \frac{ - \pi | \lambda |^2 \xi |^2}{2 ( i - i t)}}$
we see that 
$$
\hat{I}^\delta_{\lambda} ( x ) \, = \,
\displaystyle{\int^{\infty}_{- \infty}} \,
e^{ - 2 \pi it} \,
e^{- \pi \delta t^2} \,
e^{- 2 \pi \frac{| x|^2}{\lambda^2} ( 1 - i t)} \: dt \, ,
$$
which in turn is $e^{- 2 \pi | x |^2/\lambda^2} \,
\delta^{- 1/2} \,
e^{- \pi ( 1 - |x|^2/\lambda^2)/\delta} \, .  $
Inserting this in (7.2), and letting $\delta \rightarrow 0$, 
we obtain 
$$
\displaystyle{\int_{\Bbb{R}^d}} \,
\widehat{\varphi} ( \xi ) \,
I_\lambda ( \xi ) \, d \xi \, = \, 
c_d \,
\displaystyle{\int_{\Bbb{R}^d}} \,
\varphi ( x ) \, d \sigma_\lambda ( x )\, ,
$$
and thus $I_\lambda ( \xi ) \, = \,
c_d \, d \widehat{\sigma}_\lambda ( \xi )$, as was to be proved.

Note that
$$
c_d \, = \, I_\lambda ( 0 ) \, = \,
e^{2 \pi} \,
\displaystyle{\int^{\infty}_{- \infty}} \,
e^{- 2 \pi i t} \,
\frac{d t}{(2 ( 1 - it ))^{d/2}} \, = \,
\frac{\pi^{d/2}}{\Gamma ( d / 2 )} \, .
$$

\section{Counter-examples}

Since we shall be dealing with all $d \geq 2$, we return 
to the original definition of the averages $A_\lambda$,
$$
A_\lambda ( f ) ( x ) \, = \,
\frac{1}{N ( \lambda )} \, \cdot \,
\displaystyle{\sum_{| m | = \lambda}} \,
f ( n - m ) \, .  $$
Let us take $f$ to be the unit mass at the origin; i.e.
$f ( 0 ) = 1$, and $f ( n ) = 0$, if
$n \in {\Bbb{Z}^d}$, $n \ne 0$.  Then clearly
$f \in \ell^p ( {\Bbb{Z}^d)}$, for every $p$.  Next, we observe that 
$$
A_\lambda ( f ) ( n ) \, = \,
1/_{N ( \lambda )} \, = \,
1/_{N ( | n |)} \, , \ \ {\rm if} \ \ 
| n | \, = \, \lambda \, .
$$
Hence,  
\setcounter{equation}{0}
\begin{equation}
A_\star \, ( f ) ( n ) \, = \,
\displaystyle{\sup_{\lambda}} \,
A_\lambda ( f ) ( n ) \, \geq \, 1/_{N ( | n | )} \, .
\end{equation}
(Recall that $N ( \lambda ) = $ number of $m \in {\Bbb{Z}^d}$, 
so that $| m | = \lambda \, ;$ i.e. 
$N ( \lambda ) = r_d ( \lambda^2)$.)

Consider now the situation when $d \geq 5$.  Then as we have pointed
out, $N ( \lambda ) \approx \lambda^{d - 2}$, and so
$A_\star ( f ) ( n ) \geq c | n |^{- d + 2}$.  But the latter
function belongs to $\ell^p ( {\Bbb{Z}^d})$ only when 
$p > \frac{d}{d - 2}$,
and so the necessity of that condition is proved.

Next assume $d \leq 4$.
We shall use the fact that $r_4 ( 2^{2k} ) = 24$, for every\break $k \geq 1$.
This follows from the Jacobi formula which states that
$r_4 ( m ) = 8 \cdot \sigma_1^\ast ( m )$, where $\sigma_1^\ast ( m )$ is
the sum of the divisors of $m$ which are not
divisible by $4$. (See [HW, Chap.\ 20].)
Thus
$$
r_2 ( 2^{2k} ) \, \leq \,
r_3 ( 2^{2k} ) \, \leq \,
r_4 ( 2^{2k}) \, = \, 24 \, . 
$$
Now for each $d$, $d \leq 4$, we then have 
$N ( \lambda ) \leq 24$, if $\lambda = 2^k$.  And so for
$n \in {\Bbb{Z}^d}$ with $n = ( 2^k , 0 , \ldots )$, we see
that $A_\star ( f ) ( n ) \, \geq \,
1/24$, by (8.1).  Because this happens for infinitely many $n$, we have
$A_\star ( f ) \notin \ell^p$, for any $p < \infty$,
and so the necessity of the condition $d \geq 5$ is established.

\AuthorRefNames [SHW]

\end{document}